\documentclass[]{interact}

\usepackage{epstopdf}
\usepackage[caption=false]{subfig}

\usepackage[numbers,sort&compress]{natbib}
\bibpunct[, ]{[}{]}{,}{n}{,}{,}

\theoremstyle{plain}

\theoremstyle{definition}

\theoremstyle{remark}

\usepackage{listingsutf8,newtxtt,xcolor,graphicx,xurl}
\makeatletter
\def\lst@outputspace{{\ifx\lst@bkgcolor\empty\color{white}\else\lst@bkgcolor\fi\lst@visiblespace}}
\makeatother

\begin{document}

\articletype{ARTICLE TEMPLATE}

\title{Towards understanding the central limit theorem by learning Python basics}

\author{
\name{Zolt\'an Kov\'acs\textsuperscript{a}\thanks{CONTACT Z. Kov\'acs. Email: zoltan@geogebra.org}
and Alexander Thaller\textsuperscript{b}}
\affil{\textsuperscript{a}
The Private University College of Education of the Diocese of Linz,
Salesianumweg 3, Linz, Austria;
\textsuperscript{b}Linz School of Education,
Altenberger Stra\ss e 54, Linz, Austria}
}

\maketitle

\begin{abstract}

We report on a first experiment about an email based course that connects learning
Python basics and introductory probability theory. In the experiment 7 short sequences
of homework were sent out to prospective mathematics teachers who did not have
any programming background formerly, but already had some minor knowledge on probability theory.
The experiment was about to decide if learning basics of programming can promote
understanding main concepts of probability theory.
\end{abstract}


\section{Introduction}

We, the authors, have been spending a reasonable time with private tutoring in many fields
of mathematics---we are teachers. And, we agree on that probability theory is maybe
the most difficult field of mathematics---if it is about \textit{understanding}. For us, rough explanation
of statistical relationships seems to be easier than explaining the delicate issues
of properties of binomial coefficients, infinite sums and convergence. Indeed, in some
sense, statistics is much easier to explain---one needs to run virtual experiments
in a web browser
and this is well supported for several years, among others, in the online version of GeoGebra\footnote{%
See, for example, Steve Phelps' online applets
\url{https://www.geogebra.org/m/bxytp4hq},
\url{https://www.geogebra.org/m/qsjbqjfe} and
\url{https://www.geogebra.org/m/bu5sxbn2},
that visualize experiments with ice cream purchases, a broken stick or a fair coin, respectively;
or check out Andreas Lindner's collections at \url{https://www.geogebra.org/m/AytaSakt}
and \url{https://www.geogebra.org/m/qXPnyCTZ}.}.

Connecting pure mathematics and real life or virtual experiments, can be, to our experience,
extremely hard. Probability is an important item of the curriculum at secondary level,
thus it is also an important item in the teacher training. But, to our experience,
many prospective teachers have difficulties \textit{even with the basics} beyond introductory level.
Our conjecture is that the main problem lies in students' verification if theory indeed meets practice,
so, at the end of the day, a very uncertain knowledge can be observed in the students' perception. As a consequence,
in most cases, not just yesterday's and today's mathematics teachers cannot connect
theory and practice, but their students either. Roughly speaking: everybody talks about probability and
statistics, but just a few have \textit{real understanding} behind the scenes.

In our paper we propose a radical way, or at least an extension of the classic way, how
basics of elementary probability theory should be taught. It is \textit{programming}.
Our conjecture is that understanding will be significantly improved if the main concepts
of probability theory are supported with simple computer programs. The output for their flexible input
should be immediately checked, eventually in a web browser, and by changing the input parameters
to a higher number, the results can be quickly generalized. We emphasize that computers
can take on the high number of computations from human---and this is what we exactly want.
A quite simple \textit{sample space} can contain \textit{a lot} of elements, even millions or more, and it cannot
be expected that non-experts can count their elements without deeper knowledge. In fact,
experts in probability are usually experts in combinatorics as well. To reach a reasonable
level of understanding of the main concepts in probability, by using just paper and pencil,
one needs to have a very strong and \textit{safe} background in combinatorics.

In this paper we demonstrate our proposal by explaining an email based course, sent out
to 4 prospective mathematics teachers, allowing them an arbitrary time to solve the problems being set.
The students claimed that they did not have any deeper knowledge in programming formerly, but all of
them already studied probability theory at university level. Their knowledge, however, was
not yet checked via examination. The course consisted of 7 sheets of homework assignments
and took place between November 2020 and January 2021, at the Johannes Kepler University of Linz, Austria.

\section{Preparations}
\label{preps}

First author hold a series of supplementary lectures on probability theory two years ago, where second author
took part as a student and illustrated several homework assignments by using Python programs.
Some of his illustrations were a bit complicated for beginners, so the idea came that
simplified versions of such programs could be helpful for many students.
But, instead of giving a problem in probability theory and asking for a solution, in this
project we reversed the logic: We revealed the Python code in advance, and wanted to know
what problem will be solved with its help.

\subsection{Why Python?}

Why did we choose Python as the programming language for this course?
There are a couple of reasons we would like to note. First of all: It is
free, independent of the used operating system, has good documentation
and great tutorials. In addition, Python is currently a very popular
programming language: according to
\url{https://madnight.github.io/githut/#/pull_requests/2020/4}, Python
was the second most favoured language on GitHub in the 4th quarter of 2020;
moreover, ZDNet's recent article \cite{LiamTung2020} highlights that Python
was the second best liked language in November 2020.
Many programming courses teach it as a first
programming language because of its ease of use and short expressive
code.

Also, we think that Python is a great tool to express mathematical
concepts in a concise and direct way---several recent books give
examples on a possible start, see e.g.~\cite{Fuhrer2016}, \cite{Mehta2015} and \cite{Stewart2014}.
In fact, in our course students
needed to learn classic programming concepts like conditions and loops
together with probability, so we needed a programming language with a
steep learning curve. We think Python can be seen as a great tool that
supports students with the understanding of mathematical probability
concepts without the need to learn too much programming.

On the other hand, in the field of data science
Python is currently the mostly used programming language. Especially in
machine learning, which could be seen as an advanced application of
probability, Python is very well liked. We think that students who believe
that the exercises they solve are not artificially constructed for the
course, but have some connection to the real world, will be more
motivated. As a teacher we can give students the impression that if they
use Python in their introductory course on probability, they take a first
step into the prestigious field of data science. This may increase the
motivation of students and as a result they will achieve a higher level.

\subsection{The programs}

\begin{figure}
\begin{center}\includegraphics[scale=0.4]{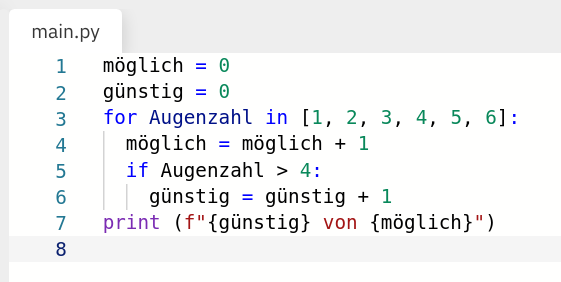}
\caption{An explanatory screenshot to help the students' first steps}
\label{main-py-window}
\end{center}
\end{figure}

The first example of the introduced Python programs
is shown in Fig.~\ref{main-py-window}. We kept the original version to highlight
the possibility to write special characters in the Python code (we used Python version 3),
including accented German letters. In the first sheet of problems each student had to try
the code on his or her own on the website \url{https://repl.it/languages/python3},
by paying attention on the correct indentation as well---the students, however, did not have
to type the code but just copy it from the email into the programming window.
We also instructed all students with an
explanatory screenshot how the program should be started and the output can be obtained.

In this introductory sheet the first question was to make the program work by running it,
and then explain its meaning line by line.
To give a somewhat creative assignment we also asked how the output $2/6$ can be changed to
the $3/6$. The outcome of this puzzle will be analyzed later in Section \ref{problem-sheets}.

It was expected that the students solve the sheet problems on their own, without
any communication. However, external help like using Internet was allowed for them.
The students were informed that about 6--8 sheets should be solved in the project,
and the required time for each should not exceed one hour.

In fact, we created a couple of other programs in advance, but always waited for
all solutions to be turned in by the students before finalizing the next sheet of problems.
We \textit{did not} give any feedback on their solutions. Instead, if we felt that
some difficulties arose, then we tried to cover them in the next sheet of problems.
The reason behind this was to test if it is possible to create a set of problems
\textit{in advance} for a possible \textit{static} learning path. Even if our experiment
was not fully static, for a future project now we already have a static set of materials.

Here we emphasize that easy access to an online version of a Python interpreter
played a significant role in our project. We did not want to ask the students
to install anything on their computers. Instead, we tried to minimize the gap
between theory and practice. In fact, \url{repl.it} is not the only website
allowing the user to try his or her code freely---we were just
satisfied with its simple user interface. Later, to obtain graphical output in an easy way, we also
introduced how Sage worksheets at \url{cocalc.com} can be accessed---in this latter case,
however, accented letters did not look elegant.

Importantly, we simplified and minimized the code we used during the email course.
This pre-process required several discussions among the authors to avoid
the technical difficulties of the programming language, but still allow
obtaining simple outputs. We had to balance between readability, simplicity,
elegance and speed. We learned that Python itself could be improved to
provide a simpler way to construct iterations than the usual \texttt{for i in range(...)} loop.
At the end of the day we were still satisfied with all of the code and the students'
solutions confirm that our choices were acceptable.

As an example of a successful balance, we mention that we heavily used Python's support
for Unicode variable names.
We denoted the sample space by $\Omega$ and a single outcome by $\omega$ since
these notations were already well-known for the students. On the other hand,
typing these Greek characters is not straightforward: one may need to copy-paste
them since they are not available on a German keyboard.

\section{The problem sheets}
\label{problem-sheets}

\subsection*{Sheet 1}
The first sheet has already been mentioned in Section \ref{preps}.
The expected explanation of the Python program was solved well by the students.
All of them managed to write a correct answer to the creative question, namely
\begin{lstlisting}[language=Python,firstnumber=5]
  if Augenzahl > 3:
\end{lstlisting}
One student mentioned that the same output can be obtained with this change:
\begin{lstlisting}[language=Python,firstnumber=2]
günstig = 1
\end{lstlisting}
when line 5 is kept in its original form. Of course, this answer makes sense
only from the programmer's point of view, and it has less to do with probability.

After receiving the results for the first sheet we had the impression that
the students got familiar with the basics of the required technology.

\subsection*{Sheet 2}

The second sheet consisted of the following four questions:

\begin{enumerate}
\item Consider the following program:
\begin{lstlisting}[language=Python]
total = 0
desired = 0
for outcome_die_1 in [1,2,3,4,5,6]:
  for outcome_die_2 in [1,2,3,4,5,6]:
    total = total + 1
    if outcome_die_1 + outcome_die_2 == 10:
      desired = desired + 1
print(f"{desired} of {total}")
print(f"p = {desired/total}")
\end{lstlisting}

\item Explain its meaning line by line, in one sentence for each.

\item Give the probability space where the possible outcomes are counted.
Which event has been identified here as ``desired''?

\item Modify the program to get a mechanical solution for the following question:

\begin{quote}
One fair die is rolled three times. Compute the probability for the event that the outcomes
correspond to a monotone increasing sequence.
\end{quote}

\end{enumerate}

We asked the students to try to solve the sheet problems without using any external help.
In case of lack of success
they were allowed to use a website that explains how to use if-then-else conditions in Python.

All students faced the quoted question during the introductory course in probability, but that time the solution
was expected to follow a purely pencil-and-paper method.

\subsubsection*{The results}

This sheet was much more difficult than the first one. One student reported that about an hour was required to solve it.
We learned that one student has a misconception by mixing relative frequency and probability---it was
a kind of surprise for this student that the program computes \textit{probability}, even if
there is no randomness built in. (In this case the students' misconception was corrected
by an explanation given in an email answer. Here we made an exception and forced
communication to clarify this basic concept.)

In general we were satisfied with the students' answers. All of them learned how three loops can
be nested to iterate on $\Omega=\{1,2,\ldots,6\}^3$. Multiple ways were found to express monotonicity,
including both the construct with an \texttt{and} clause and the short form
\texttt{outcome\_die\_1 <= outcome\_die\_2 <= outcome\_die\_3} which is one of Python's strengths.

In our experience it turned out that with little work a definitely large sample space
can be studied. In fact, manual check of $\{1,2,\ldots,6\}^3$ (or, equivalently, deriving the sum
$1+3+6+10+15+21$) is much more time consuming than writing a suitable program.

\subsection*{Sheet 3}

Fig.~\ref{3.py} shows the main code that was announced in the third sheet.

\begin{figure}
\begin{lstlisting}[language=Python]
import itertools
# 1. The basic set:
outcomes = {1, 2, 3, 4, 5, 6}
# itertools.product creates a Cartesian product:
Ω = set(itertools.product(outcomes, outcomes))
print(Ω)
# Note that Ω as a set and it is unsorted.
# Instead of "set" we can use a "list":
# in that case the data will be sorted.

# 2. Count the amount of desired and total outcomes:
desired = 0
total = len(Ω)
for outcome_die_1, outcome_die_2 in Ω:
  if outcome_die_1 + outcome_die_2 == 10:
    desired = desired + 1

# 3. Output:
print(f"{desired} of {total}")
print(f"p = {desired/total}")
\end{lstlisting}
\caption{The third program (translated into English)}
\label{3.py}
\end{figure}

With this program we introduced some handy techniques via \texttt{itertools.product}
to define Cartesian powers. We pointed to
the second sheet and asked the students to rewrite their programs by using the new techniques.
We also asked them to display all desired outcomes on the screen (with the \texttt{print} command)
in a sorted fashion (via the \texttt{list} datatype).

We had positive feedback---they managed to understand the point. A student successfully generalized the 
solution for the puzzle of the second sheet for 5 dice, even if it was not directly asked.

\subsection*{Sheet 4}

The fourth sheet consisted of the following four questions:

\begin{enumerate}
\item Consider the following program:
\begin{lstlisting}[language=Python]
import itertools

def X(ω):
  """Students could edit this function on their own.
  X is a rendom variable, that is, X: Ω → ℝ.
  A real number will be mapped to each outcome."""
  return ω[0]+ω[1]

outcomes = {1, 2, 3, 4, 5, 6}
Ω = list(itertools.product(outcomes, outcomes))
# Event E
E = []
for ω in Ω:
  if X(ω) == 10:
    E.append(ω)
print(f"{len(Ω)} total outcomes")
print(f"{len(E)} desired outcomes")
print(f"p = {len(E)/len(Ω)}")
\end{lstlisting}

\begin{enumerate}
\item Check if the program works properly. Save a screenshot for an evidence.
\item Explain the difference between this program and the one in sheet 3.
\end{enumerate}

\item Modify the program on line 7 and 14 in order to solve this problem:
\begin{quote}
Alice and Bob roll a die four times consecutively.
If the sum of the outcomes is between 7 and 14, Alice wins, otherwise Bob.
For which player is the game advantageous?
\end{quote}

\item Change the rules of the game above to get a fair play.
Verify the changed rules by writing a suitable program.

\item Assume that lines 12--15 are replaced by the following ones:
\begin{lstlisting}[language=Python,numbers=none]
E = [ω for ω in Ω if X(ω) == 10]
\end{lstlisting}
Check that the same result can be obtained as before.
Also, explain this short form by using a mathematical formula like a definition for a set.
\end{enumerate}

In this sheet we focused on the concept of random variables. To highlight their importance
we defined a Python function \texttt{X($\omega$)}.

\subsubsection*{The results}

The students answered the questions generally quite well. Two of them remarked that
the problem setting in question 2 was incomplete since line 10 has to be changed as well.

It turned out that the ``copy-paste'' way of copying a program was technically not always
straightforward---sometimes indentation was transferred improperly and some manual edit
was required on the student's side.

\begin{figure}
\begin{center}\includegraphics[width=1.0\textwidth]{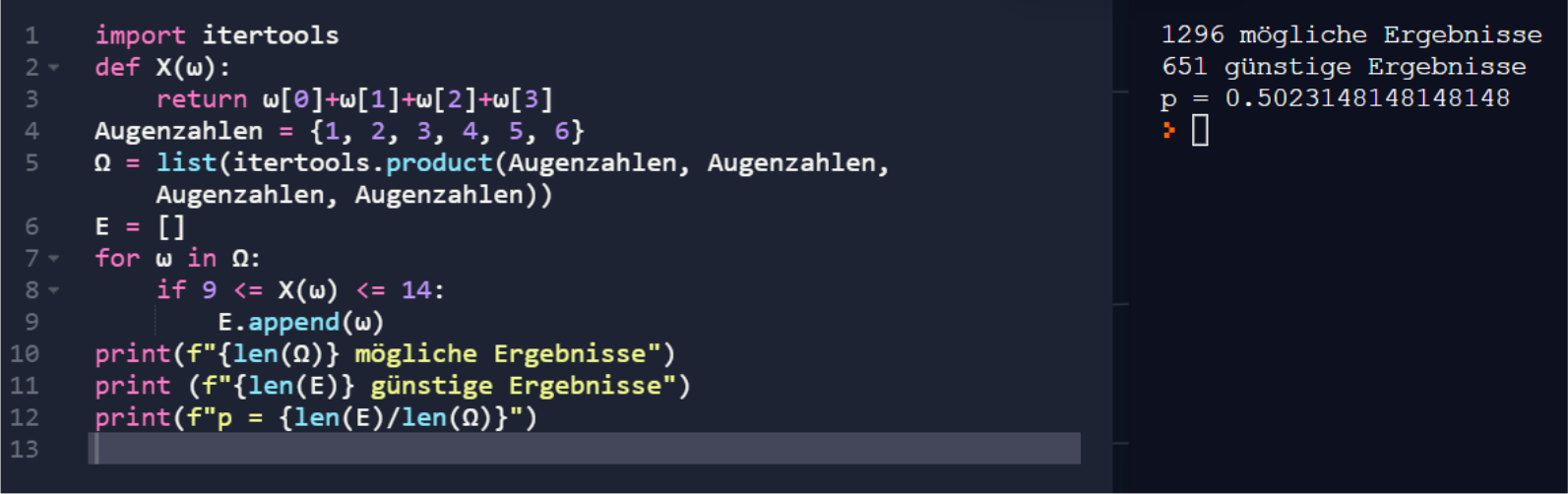}
\caption{A student's solution for question 3 in sheet 4}
\label{42}
\end{center}
\end{figure}

Fig.~\ref{42} shows a solution by a student for question 3. It is clear that this
student successfully combined previous programming knowledge, theory and practise. Similarly,
two other students found the same idea to change the limits to have $7\leq X(\omega)\leq14$,
and one of them tried several attempts to find better limits---without success.
(Note that for an even number of dice there is no trivial solution by halving the
sample space symmetrically in the middle. So, in some sense, question 3 is adequate to motivate
experimenting.)

The last question was solved perfectly by 3 students, giving the answer
$$E=\{\omega\;\mid\;\omega\in\Omega\;\land\;X(\omega)=10\}$$
where $X(\omega)=\omega_1+\omega_2$. Here we concluded that the technical
fact that Python's indexing starts from $0$, can be confusing for beginners,
since in mathematics the usual indexing starts from $1$.

\subsection*{Sheet 5}

Fig.~\ref{5a.py} shows the first program in the fifth sheet. The students' assignment was
to find a text problem that is solved by the program.

\begin{figure}
\begin{lstlisting}[language=Python]
import itertools
Ω = list(itertools.product({"h", "t"}, repeat=5))
print(f"Ω={Ω}")
E = [ω for ω in Ω if ω.count("h") == 2]
print(f"E={E}")
print(f"{len(E)} of {len(Ω)}")
print(f"p = {len(E)/len(Ω)}")
\end{lstlisting}
\caption{The first program in the fifth sheet (translated into English)}
\label{5a.py}
\end{figure}

To find a connection with binomial coefficients, the students had to explain why the result
equals to $\binom{5}{2}/2^5$. A technical question helped in improving the look of the output,
by suggesting the addition of this line between lines 2 and 3:
\begin{lstlisting}[language=Python,numbers=none]
Ω = {"".join(outcome) for outcome in Ω}
\end{lstlisting}

Then a second program was introduced. In fact, it was just a slight modification of the first program,
but, for simplicity, the whole code was announced (see Fig.~\ref{5b.py}).
\begin{figure}
\begin{lstlisting}[language=Python]
import itertools
Ω = set(itertools.product({"h", "t"}, repeat=5))
expected_value = 0
for ω in Ω:
  expected_value = expected_value + ω.count("h")
expected_value = expected_value / len(Ω)
print(f"expected value = {expected_value}")
\end{lstlisting}
\caption{The second program in the fifth sheet (translated into English)}
\label{5b.py}
\end{figure}
Here the assignment was to give a random variable $X$ whose expected value is computed in the program.
Another assignment was to explain why the expected value $2.5$ is printed. Then, this result
had to be generalized by changing the code \texttt{repeat=5} to other values. Finally,
also the variance had to be computed by adding some more lines to the code.

\subsubsection*{The results}

The first program was well identified by all students---a fair coin was tossed 5 times,
and the probability of the appearance of exactly two heads was computed.
Applying the technical hint was unsuccessful for two students (for one the closing brace was missing,
for another one the issue was unclear).

For the second program, one student added the line \texttt{print(len($\Omega$))} to better understand the
code. (In fact, the function \texttt{len} was already introduced in the third sheet.) Instead
of mentioning symmetry, the same student argued by citing the formula $E=n\cdot p$ for the
binomial distribution where $n$ denotes the number of independent experiments and $p$ the probability
of success (here $0.5$). The same student computed the variance by multiplying $E$ with $0.5$---this
is correct for the given assignment, but here no programming was used, just theory.

\begin{figure}
\begin{center}\includegraphics[width=1.0\textwidth]{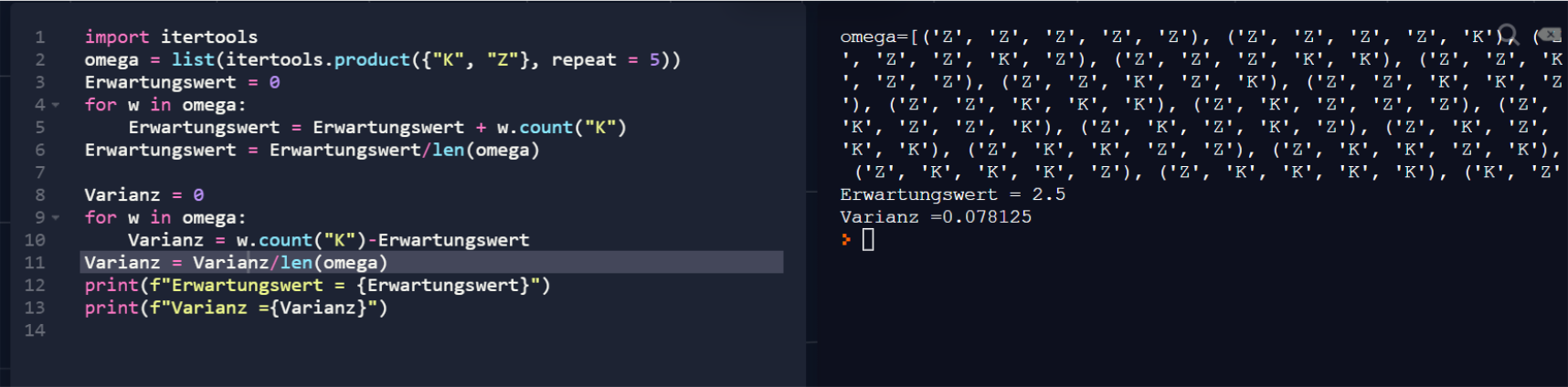}
\caption{A student's incorrect program code to compute variance in sheet 5}
\label{52}
\end{center}
\end{figure}

The other students explained the result $2.5$ by mentioning symmetry, and the same argument was
used for the generalization. But no answer was found to compute the variance, only for one student,
and in that case it was incorrect (see Fig.~\ref{52})---clearly, on line 10 something like
\begin{lstlisting}[language=Python,firstnumber=10]
  Varianz = Varianz + (w.count("K")-Erwartungswert)**2
\end{lstlisting}
should have been written. Actually, the notation \texttt{**2} could be something unfamiliar
for those who never did any programming before. Anyway, the result $0.078125$ should have been
suspicious. It is also remarkable that this student preferred the Latin alphabet to Greek characters.

\subsection*{Sheet 6}

We identified the assignment to compute the variance in the fifth sheet somewhat unsuccessful.
So we continued by sharing two correct solutions to compute the variance (see Fig.~\ref{6a.py}
and \ref{6b.py}).

\begin{figure}
\begin{lstlisting}[language=Python]
import itertools
Ω = set(itertools.product({"h", "t"}, repeat=5))
E = 0
for ω in Ω:
  E = E + ω.count("h")
E = E / len(Ω)
V = 0
for ω in Ω:
  V = V + (ω.count("h")-E)**2
V = V / len(Ω)
print(f"expected value = {E}")
print(f"variance = {V}")
\end{lstlisting}
\caption{A correct solution to compute the variance (translated into English)}
\label{6a.py}
\end{figure}

\begin{figure}
\begin{lstlisting}[language=Python]
import itertools
Ω = set(itertools.product({"h", "t"}, repeat=5))
E = 0
V = 0
for ω in Ω:
  x = ω.count("h")
  E = E + x
  V = V + x*x
E = E / len(Ω)
V = V / len(Ω) - E*E
print(f"expected value = {E}")
print(f"variance = {V}")
\end{lstlisting}
\caption{Another correct solution to compute the variance (translated into English)}
\label{6b.py}
\end{figure}

The students had to compare these two program codes, check if their outputs are the same
(even for 6 or 7 dice), and find the differences between the way of programming.
We actually wanted to add another connection with theory, namely, to verify if the
formula
\begin{equation}
V(X)=E((X-E(X))^2)=E(X^2)-E(X)^2
\label{variance}
\end{equation}
(also known as Steiner translation theorem in Germany and König-Huygens formula in the French literature)
can be reflected in the two solutions.

As another activity, a third program was given (see Fig.~\ref{6c.py}) to prepare the last
topic, namely, some steps towards the central limit theorem (CLT, \cite{Polya1920}).
\begin{figure}
\begin{lstlisting}[language=Python]
import itertools
Ω = set(itertools.product({"h", "t"}, repeat=5))
possible_values = [0, 1, 2, 3, 4, 5]
# In the beginning each value of the random variable
# has a frequency of 0.
heads = [0]*len(possible_values)
for ω in Ω:
  heads[ω.count("h")] += 1
probabilities = [h/len(Ω) for h in heads]
print(list(zip(possible_values, probabilities)))
\end{lstlisting}
\caption{A third program in the sixth sheet (translated into English)}
\label{6c.py}
\end{figure}
Here we provide the questions to this third program:
\begin{enumerate}
\item What computation is performed by this program?
\item Consider following modifications:
\begin{lstlisting}[language=Python,firstnumber=2]
Ω = set(itertools.product({1, 2, 3, 4, 5, 6}, repeat=4))
possible_values = range(30)
\end{lstlisting}
\begin{lstlisting}[language=Python,firstnumber=8]
  heads[sum(ω)] += 1
\end{lstlisting}
Explain the output of the modified program.
\item We are about to plot the obtained data. So we change line 10 as follows:
\begin{lstlisting}[language=Python,firstnumber=10]
for w in zip(possible_values, probabilities):
  print(f"{w[0]},{w[1]}")
\end{lstlisting}
Copy the output of this program in a spreadsheet software and create an XY diagram
to get some graphical output.
\item Extend the program with these two additional lines:
\begin{lstlisting}[language=Python,firstnumber=12]
distribution = [sum(probabilities[:h]) for h in possible_values]
print(list(zip(possible_values, distribution)))
\end{lstlisting}
Explain the result.
\end{enumerate}
In fact, the variable naming was somewhat incorrect here (because \texttt{heads} 
actually refers to the sums of the outcomes), but we wanted to make very little
changes on the code.

Among the authors there was a quite a long discussion on getting the output in
the best fashion. The Python construction \texttt{list(zip(...))} is a kind
of compromise between beautiful output with much code and simple output with less code,
however, this is something non-trivial in Python.
Finally a graphical output was also expected---students should be familiar with
popular spreadsheet software and visualize data in an easier way than programming.

\subsubsection*{The results}

For the first part, all solutions were correct: three of them mentioned that
for the second code only one loop is required, so it is somewhat more effective.
In fact all of the turned-in solutions lacked of the wow feeling that we used
the formula (\ref{variance}). One student reported
experiments to learn if \texttt{**} is a typo or not---finally the correct interpretation
was found: it is \textit{the power sign}.

For the second part, the students successfully identified the purpose of the third program
(that is, to find all probabilities for the different number of heads for 5 dice), and of the changed one
(to find all probabilities of the sums of the outcomes for 4 dice). Both the density
function and the cumulative distribution function (CDF) were obtained correctly, and for the former,
a graphical output was always attached. (One student created a graph for the latter as well.)
We also learned that, in the case of three students, the density function of the
\textit{discrete} binomial distribution
was incorrectly identified, namely, as a \textit{continuous} normal distribution.

We also noted that the modification \texttt{possible\_values = range(30)} aroused interest
among the students,
namely, to find the meaning of the number 30. (In fact, only the sums $4,5,\ldots,24$
had a meaning, but for simplicity we wanted to avoid any manipulation of the indexes of
the arrays.)

\subsection*{Sheet 7}

In the final sheet we wanted to point towards CLT, and thus we gave two programs:
they can be seen in Fig.~\ref{7a.py} and \ref{7b.py}).

\begin{figure}
\begin{lstlisting}[language=Python]
import itertools
max_dice = 4
for d in range(1, max_dice + 1):
  Ω = set(itertools.product({1, 2, 3, 4, 5, 6}, repeat=d))
  sums = [0] * (max_dice * 6 + 1)
  for ω in Ω:
    sums[sum(ω)] += 1
  probabilities = [h/len(Ω) for h in sums]
  print(probabilities)
\end{lstlisting}
\caption{First program in the seventh sheet (translated into English)}
\label{7a.py}
\end{figure}

\begin{figure}
\begin{lstlisting}[language=Python]
import itertools
import matplotlib.pyplot as p
max_dice = 4
for d in range(1, max_dice + 1):
  Ω = set(itertools.product({1, 2, 3, 4, 5, 6}, repeat=d))
  sums = [0] * (max_dice * 6 + 1)
  for ω in Ω:
    sums[sum(ω)] += 1
  probabilities = [h / len(Ω) for h in sums]
  p.plot(range(1, max_dice * 6 + 1), probabilities[1:])
p.show()
\end{lstlisting}
\caption{Second program in the seventh sheet (translated into English)}
\label{7b.py}
\end{figure}

\begin{figure}
\centering
\subfloat[An XY diagram via \texttt{pyplot.plot}]{%
\resizebox*{6cm}{!}{\includegraphics{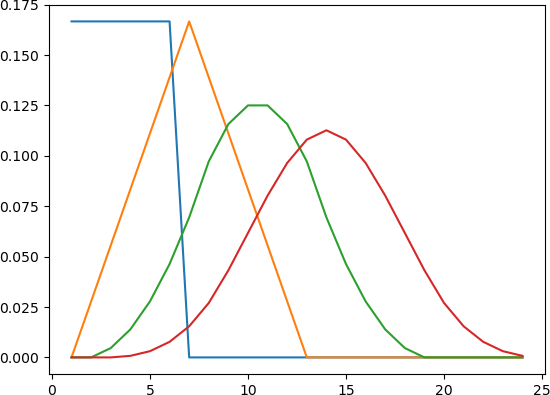}}\label{72a}}\hspace{5pt}
\subfloat[A bar chart via \texttt{pyplot.bar}]{%
\resizebox*{6cm}{!}{\includegraphics{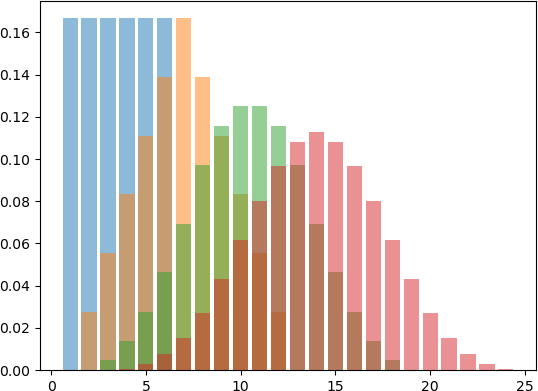}}\label{72b}}
\caption{Graphical outputs of the second program in the seventh sheet} \label{sample-figure}
\end{figure}

In fact, these two programs are almost the same, just the outputs differ: the second
one results in a nice diagram via \texttt{matplotlib} (see Fig.~\ref{72a}). We provided an alternative
way to plot a vertical bar chart with the following suggestion:
\begin{lstlisting}[language=Python,firstnumber=10]
  p.bar(range(1, max_dice * 6 + 1), probabilities[1:], alpha=0.5)
\end{lstlisting}
This is mathematically more correct because the discrete values are highlighted
(see Fig.~\ref{72b}).
On the other hand, to point towards CLT, also a continuous diagram makes sense.

Technically, we asked the students to use
Sage worksheets at \url{cocalc.com} since it provides a simple interface to obtain
the graphical result.

In this closing activity we requested to obtain as many mathematical properties of the bar chart
as possible. We asked the students to generalize their experiences. We wanted them to
focus on using notions which are also familiar in a secondary school context. We clarified
that notions like sequence, limit, function, density function, expected value, variance,
standard deviation and normal distribution are all allowed.

\subsubsection*{The results}

The students managed to visualize both diagrams and most of them concluded that the bar chart is
more useful to study the exact values of the functions. A student, however, preferred the XY diagram
that may provide more precise values.
They managed to read off the following relationships:
\begin{enumerate}
\item The definite integral of the density function is $1$ (between $0$ and $25$).
\item The maximum of the density function is reached at the expected value. There is a clear
relationship between the expected value and the number of rolls.
\end{enumerate}

According to CLT, the following statements were claimed (here and from now on we use $d$
for the number of dice):
\begin{itemize}
\item ``The more rolls you perform, the clearer the normal distribution of the [density] function becomes.''
\item For $d=3$ and $d=4$ ``it is clear that it is about the normal distribution of the sums,
because of its beautiful bell shaped curve.''
\item For $d=3$ ``we already obtain an approximate Gaussian bell curve (normal distribution).''
\item For $d=4$ ``we can already see a beautiful normal distribution with the expected value 14.''
\end{itemize}

Unfortunately, the notion ``limit'' was mixed with the notion ``maximum'' in some sense, so two students did not
observe at all how the situation changes if $d\to\infty$. One student studied each density function
in the case $h\to\infty$ (cf.~line 9 in Fig.~\ref{7b.py}), to obtain its limit, which is, of course, $0$.
(Also the same student dealt with the case $h\to-\infty$.) Finally, this student used 
the $68.27\%$ rule (by using external information from a different course on
probability theory) to get the standard deviation in the case $d=4$,
under the assumption that $i$ is an inflection point of the density function:
\begin{equation}
\Phi\left(\frac{i-E(X)}{\sqrt{V(X)}}\right)=0.6827,
\end{equation}
where $\Phi$ denotes the CDF of the standard normal distribution, and the inflection point is guessed
to be at $i=18$. With this approximation formula the value $\sqrt{V(X)}\approx8.33$ can be obtained.

\section{Discussion}

We did not have any preliminary goals at the beginning of this project, only to connect theory and
practice via programming. After evaluating the solutions we have the feeling that all students
have a safer approach to problem solving in probability theory, and they have no fear to give
a computer program a chance to do experiments. That is, even if they need to modify the program code,
studying via programming experiments is a feasible way.

Our project confirmed that Python is a very useful and practical tool to collect exact data
to study probability theory---even if there can appear some minor inconveniences.
As one student pointed out, by using Sage, symbolic computation
supports even more precise observations---for example, the first program in the last sheet
(Fig.~\ref{7a.py}) was accidentally launched in Sage as well and all results were shown
as vulgar fractions.

On the other hand, connecting theory and the result of a program was not always straightforward.
Students identified the CDF very well, but they were unable to recognize the CLT, even if
it was discussed during the former course---but maybe without a result of deeper understanding.
On the other hand, the CLT was not directly asked in the seventh sheet,
so the problem setting was maybe not clear enough.

At this point we need to give some description on which form of CLT should be expected to discover by the students.
In fact, CLT has several versions: both formal and informal variants. Roughly speaking, CLT can be summarized
in the following form (see \cite{Wikipedia-CLT}):
\begin{quote}
When independent random variables are added, their properly normalized
sum tends toward a normal distribution (informally a bell curve) even if
the original variables themselves are not normally distributed.
\end{quote}
A formal definition given in \cite[p.~241]{Montgomery2014} is as follows:
\begin{quote}
If $X_{1},X_{2},\ldots,X_{n}$
are $n$ random samples drawn from a population
with overall mean $\mu$ and finite variance $\sigma^2$,
and if ${\bar {X}}_{n}$ is the sample mean, the
limiting form of the distribution, $$Z=\lim _{n\to \infty}{\sqrt {n}}{\left({\frac {{\bar {X}}_{n}-\mu }
{\sigma}}\right)},$$ is the standard normal
distribution.
\end{quote}
In fact, we avoided going into any formal detail during our course by assuming that the provided formula
here \textit{may not promote understanding}.

\begin{figure}
\begin{center}
\includegraphics[width=0.7\textwidth]{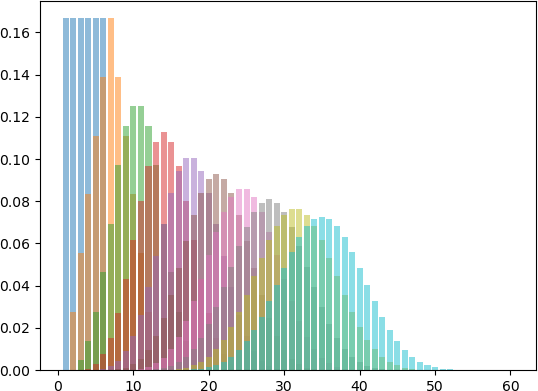}
\caption{Increasing the value of \texttt{max\_dice} to $10$}
\label{10}
\end{center}
\end{figure}

More importantly, we conjecture that making a difference between a discrete distribution and a continuous one may be
a big jump in the level of understanding for the student. Discrete distributions are \textit{never} normal,
and to speak about a Gaussian distribution one should introduce the notion of limit in a very
delicate form. Here, we should translate the random variable $X$ to have the expected value $0$ first,
and then normalize it by a function of $d$. We did not see any sign of this concept among the
students' solutions. Maybe, in a next project, we will explicitly ask the students to try to
\textit{approximate} the discrete distributions with a continuous one, by searching for the adequate
expected value and standard deviation. Also, experimenting with higher values for $d$ should also be
suggested---no student
tried to increase the value of \texttt{max\_dice}, even if technically it is possible to try
\texttt{max\_dice = 8} in a \texttt{cocalc.com} Sage session. (Fig.~\ref{10} shows the output
for \texttt{max\_dice = 10} when run natively on a recent computer, after about 1 minute
of computation. Here we warn the reader that a sample space $\Omega$ with more than 60 millions of elements
must be observed, since $|\Omega|=6^{10}=60,466,176$.)

Finally, we remark that experimenting with a non-fair die could be another step to explain
the CLT in a more general approach. A fair die, after rolled three or four times,
results in a ``too symmetrical'' binomial distribution---it can be confused with
the normal distribution. A non-symmetrical result, which is still similar enough
to the normal distribution to make a conjecture, could be eventually a better
example to promote understanding the CLT and help students in formalizing it---just
by using their experiments and own words.
(For a nice continuous example we refer to the Wikipedia page
\url{https://en.wikipedia.org/wiki/Illustration_of_the_central_limit_theorem#Illustration_of_the_continuous_case}.)

\section{Conclusion}

We reported on an experimental project that connects an email based Python course
and former knowledge of elementary concepts of probability theory.
We think our results successfully present that there is room for this kind of approach
in the education. We also think that our method should be further refined and tested.

We emphasize that this research just the first step in a longer research project.
Basically, we would like to continue with interviews with the participating
students whether their knowledge and motivation were improved during the course.
This will give us further feedback to meet the best decision towards further
refinements of the course.

All materials of our project---including the Python programs, the homework assignments and
the turned-in homework---can be found online at \url{https://github.com/kovzol/probability-python}.

\section*{Acknowledgements}

We are thankful to the students Alexander, Romana, Sandra and Sarah
for their participation in the project, and to Andreas Lindner
and two anonymous reviewers who
gave us useful feedback on the first versions of the manuscript.

\bibliographystyle{tfs}
\bibliography{article}

\begin{thebibliography}{1}
\providecommand{\MR}{\relax\unskip\space MR }
\providecommand{\url}[1]{\normalfont{#1}}
\providecommand{\urlprefix}{Available at }

\bibitem{Fuhrer2016}
C. F\"uhrer, J.E. Solem, and O. Verdier, \emph{Scientific Computing with Python
  3}, Packt Publishing, 2016.

\bibitem{Mehta2015}
H.K. Mehta, \emph{Mastering Python Scientific Computing}, Packt Publishing,
  2015.

\bibitem{Montgomery2014}
D.C. Montgomery and G.C. Runger, \emph{Applied Statistics and Probability for
  Engineers}, 6th ed., 2014.

\bibitem{Polya1920}
G. P\'olya, \emph{{\"U}ber den zentralen {G}renzwertsatz der
  {W}ahrscheinlichkeitsrechnung und das {M}omentenproblem}, Mathematische
  Zeitschrift 8 (1920), pp. 171--181.

\bibitem{Stewart2014}
J.M. Stewart, \emph{Python for scientists}, Cambridge University Press, 2014.

\bibitem{LiamTung2020}
L. Tung, \emph{Programming language {P}ython's popularity: Ahead of {J}ava for
  first time but still trailing {C}},
  \url{https://www.zdnet.com/article/programming-language-pythons-popularity-ahead-of-java-for-first-time-but-still-trailing-c/}
  (2020).

\bibitem{Wikipedia-CLT}
 {Wikipedia contributors}, \emph{Central limit theorem --- {Wikipedia}{,} the
  free encyclopedia},
  \url{https://en.wikipedia.org/w/index.php?title=Central_limit_theorem&oldid=1011860578}
  (2021). [Online; accessed 13-March-2021].

\end{thebibliography}

\end{document}